\documentclass[11pt,twoside, reqno]{article}

\usepackage{amssymb}
\usepackage{amsmath}
\usepackage{mathrsfs}
\usepackage{amsthm}
\usepackage{txfonts}

\allowdisplaybreaks

\newtheorem{theorem}{Theorem}[section]
\newtheorem{lemma}[theorem]{Lemma}
\newtheorem{corollary}[theorem]{Corollary}

\theoremstyle{definition}
\newtheorem{remark}[theorem]{Remark}
\newtheorem{definition}[theorem]{Definition}

\numberwithin{equation}{section}

\begin{document}

\arraycolsep=1pt

\title{\Large\bf On the symmetric $q$-analog on the bi-univalent
functions\\ with respect to symmetric points
\footnotetext{\hspace{-0.35cm} 2010 {\it Mathematics Subject
Classification} {30C45; 30C50.}
\endgraf{\it Key words and phrases}. Fekete-Szeg\"{o} problem, Second Hankel determinant, Bi-univalent function,
Bernardi integral operator.
\endgraf This work is supported by Natural Science Foundation of Ningxia (Grant No.\
2023AAC03001) and Natural Science Foundation of China(Grant
No.\ 12261068).}}
\author{Pinhong Long, Huili Han, Halit Orhan, Huo Tang }
\date{ }
\maketitle

\vspace{-0.6cm}

\begin{center}
\begin{minipage}{13.5cm}\small
{{\bf Abstract.} Our objective is to usher and investigate the
subclass
$\widetilde{\mathcal{S^{*}_{\sum}}}^{\eta}_{q}(\mu,\lambda;\phi)$
of the function class $\sum$ of analytic and bi-univalent
functions related with the symmetric $q$-derivative operator and
the generalized Bernardi integral operator. On the one hand,
without the generalized Bernardi integral operator we estimate the
second Hankel determinants for the reduced subclasses
$\widetilde{\mathcal{S^{*}_{\sum}}}_{q}(\lambda;\phi)$ with
respect to symmetric points. On the other hand, we also give the
corresponding results of Fekete-Szeg\"{o} functional inequalities
and the upper bounds of the coefficients $a_2$ and $a_3$ for these
subclasses.}
\end{minipage}
\end{center}

\section{Introduction\label{s1}}

\hskip\parindent Let the subclass $\mathcal{S}$ of $\mathcal{A}$
be the set of all univalent functions in the open unit disk
$\Delta=\{z\in\mathbb{C}: \mid z\mid<1\}$, where the set
$\mathcal{A}$ is the class of normalized analytic functions $f(z)$
in $\Delta$ via
\begin{equation}\label{equ1.1}
f(z)=z+\sum^{\infty}_{n=2}a_{n}z^{n}.
\end{equation}

Referring to the Koebe one-quarter theorem in \cite{Du83}, we may
derive that the inverse $f^{-1}$ of $f\in\mathcal{S}$ satisfies
\begin{equation*}
f^{-1}(f(z))=z,~~(z\in\Delta)~\text{and}~
f(f^{-1}(w))=w,~~(w\in\Delta_{\rho}),
\end{equation*}
where $\Delta_{\rho}=\{z\in\mathbb{C}: \vert z\vert<\rho\}$ with
respect to the radius $\rho\geq \frac{1}{4}$ of the image
$f(\Delta)$, in particular, $\Delta_{1}=\Delta$. As is well known,
\begin{equation}\label{equ1.2}
f^{-1}(w)=w-a_{2}w^{2}+(2a^{2}_{2}-a_{3})w^{3}-(5a^{3}_{2}-5a_{2}a_{3}+a_{4})w^{4}+\ldots.
\end{equation}
If the function $f\in\mathcal{A}$ and its inverse $f^{-1}$ are
univalent in $\Delta$, then it is called to be bi-univalent.
Furthermore, denote by $\Sigma$ the class of all bi-univalent
functions $f\in\mathcal{A}$ in $\Delta$.

The symmetric $q$-derivative operator or the symmetric
$q$-difference $\widetilde{\mathcal{D}}_{q}f(z)$ (see
\cite{SAY18}) of the function $f$ is defined by
$$
\widetilde{\mathcal{D}}_{q}f(z)=\left\{\begin{array}{ll}
\frac{f(qz)-f(q^{-1}z)}{(q-q^{-1})z},~~&(z\neq0; 0<q<1),\\
f'(0),~~&(z=0)
\end{array}\right.
$$
provided that $f'(0)$ exists, and the symmetric $q$-number
 $\widetilde{[\chi]}_{q}$ is denoted by
$$
\widetilde{[\chi]}_{q}=\left\{\begin{array}{ll}
\frac{q^{\chi}-q^{-\chi}}{q-q^{-1}}~~~&\mbox{for}~~\chi\in\mathbb{C},\\
\sum^{\chi-1}_{k=0}q^{\chi-1-2k}~~~&\mbox{for}~~\chi=n\in\mathbb{N}.
\end{array}\right.
$$
We also remark that if $q\longrightarrow1_{-}$, then
$\widetilde{\mathcal{D}}_{q}f(z)\longrightarrow f'(z)$ holds,
where $f'$ is the ordinary derivative of the function $f$. As a
result, for $f\in\mathcal{A}$ we may obtain the symmetric
$q$-derivative $\widetilde{\mathcal{D}}_{q}f(z)$ as
\begin{equation*}
\widetilde{\mathcal{D}}_{q}f(z)=1+\sum^{\infty}_{n=2}\widetilde{[n]}_{q}a_{n}z^{n-1},~~~~(0<q<1).
\end{equation*}

For $\Re\eta>-1$, introduce the generalized  Bernardi integral
operator
$\mathcal{J}^{\eta}_{q}:\mathcal{A}\longrightarrow\mathcal{A}$ determined by
\begin{equation}
\mathcal{J}^{\eta}_{q}f(z)=\frac{\widetilde{[1+\eta]}_{q}}{z^{\eta}}\int^{z}_{0}t^{\eta-1}f(t)\widetilde{d}_{q}t\quad\text{for}~z\in\Delta~\text{and}f\in\mathcal{A}.
\end{equation}
Therefore, for $f\in\mathcal{A}$, it indicates that
\begin{equation}
\mathcal{J}^{\eta}_{q}f(z)=z+\sum^{\infty}_{n=1}L^{\eta}_{q}(n)a_{n}z^{n},~~(z\in\Delta),
\end{equation}
where
\begin{equation}
L^{\eta}_{q}(n)=\frac{\widetilde{[1+\eta]}_{q}}{\widetilde{[n+\eta]}_{q}}:=L_{n}.
\end{equation}
Now we point out that when $q\rightarrow1_{-}$, it is due to the
classical Bernardi integral operator $\mathcal{J}_{\eta}$
\cite{Be69}. Clearly, Alexander \cite{Al15} integral and Libera
\cite{Li65} operators are two special versions of
$\mathcal{J}_{\eta}$ for $\eta=0$ and $\eta=1$, respectively.
Besides, recently Srivastava et al. \cite{SRASA19} consider
$(p,q)$-analog of Bernardi integral operator as well as its
$q$-analog by Noor et al. \cite{NRN17}.

Up to now, the coefficient estimate problem for every
Taylor-Maclaurin coefficient $\mid
a_{n}\mid(n\in\mathbb{N}\setminus \{1,2\})$ is probably still an
open problem. Since the works of Brannan and Taha \cite{BT86} and
Srivastava et al. \cite{SMG10}, the non-sharp estimates of first
two Taylor-Maclaurin coefficients $\mid a_{2}\mid$ and $\mid
a_{3}\mid$ for some subclasses of analytic and bi-univalent
functions of class $\sum$ were investigated by many
mathematicians; refer to \cite{De13,FA11,HO12,PN15,PN17,SEA15,SGG18,XGS12,XXS12}.

Since the pioneering work by Fekete-Szeg\"{o} \cite{FS33} about
the determination of the sharp upper bounds for the subclass of
$\mathcal{S}$, Fekete-Szeg\"{o} functional inequalities were
generalized and applied into many classes of functions; refer to
\cite{AG00, AFD14, Ko87, LTW20, MB18, OMB16, PR17}. Noonan and
Thomas \cite{NT76} ever introduced the $q$-th Hankel determinant
$H_{q}(n)$ of $f$ for $q\geq1$ and $n\geq1$, and from then on
Hankel determinant was also investigated and employed by many
mathematicians; refer to Edrei \cite{Ed40} and Pommerenke
\cite{Po67}. Moreover, regarding the second Hankel determinants
for bi-univalent functions, we refer to \cite{KAZ17, MV15, SL18}
and their references therein. In fact, the Hankel determinants
$H_{2}(1)=a_{3}-a^{2}_{2}$ and $H_{2}(2)=a_{2}a_{4}-a^{2}_{3}$ are
the classical Fekete-Szeg\"{o} functional and second Hankel
determinant, respectively.

As the same as  $(p,q)$-calculus \cite{CJ91},  the $q$-calculus
\cite{Ja10,Jq10} is a generalization of the ordinary calculus
without the limit symbol, and its related theory has been applied
into mathematical, physical and engineering fields (see
\cite{GR90}, \cite{KCh02} and \cite{RMS07}). Since Ismail et al.
\cite{IMS90} firstly adopted the $q$-derivative operator to
investigate the $q$-calculus of the class of starlike functions in
disk, there had a great deal of work in this respect; for example,
refer to \cite{RASKK20, SA16, SB17, SRASA19, Uc16}. Besides, by
involving some special functions and integral/difference operators
or increasing the complexity of function classes, many new
subclasses of analytic functions associated with $q$-calculus or
$(p, q)$-calculus were studied. Here we may refer to
\cite{AAL19,DD13, FRS17, MS17, PR11} and Srivastava et al.
\cite{SAY18} for symmetric $q$-derivative operator and
\cite{SKKA19,SRASA19}. For the multivalent functions, there have
also some related advances for \cite{ASU19,Pu12,SMAZ19}.

Motivated by these works above, in the article we define and
investigate a new subclass of the function class $\sum$ of
analytic and bi-univalent functions associated with the symmetric
$q$-derivative operator and the generalized Bernardi integral
operator. Without the generalized Bernardi integral operator, we
firstly obtain the second Hankel determinant for the
$\lambda$-quasi,bi-starlike function class with respect to
 symmetric points. Besides, we also study the corresponding
Fekete-Szeg\"{o} functional inequalities and the bound estimates
of the coefficient $a_2$ and $a_3$ for these classes or their
reduced versions.

Let $f$ and $g$ be two analytic functions in $\Delta$. If there
exists an analytic function $h$ satisfying $h(0)=0$ and $\mid
h(z)\mid<1$ for $z\in\Delta$ so that $f(z)=g(h(z))$ holds, then
$f$ is called to be subordinate to $g$, i.e. $f\prec g$. Let the
series expansion form of $\phi\in\mathcal{P}$ be expressed by
\begin{equation}\label{equ2.4}
\phi(z)=1+\sum^{\infty}_{n=1}E_{n}z^{n},~~(E_{1}>0,z\in\Delta).
\end{equation}

Next, To go on our statements, the following general subclass of analytic and
bi-univalent functions associated with the symmetric
$q$-derivative operator and the generalized Bernardi integral
operator with respect to symmetric points is introduced.

\begin{definition}\label{def1}  A function $f(z)\in\sum$ given by (\ref{equ1.1}),
is in the class
$\widetilde{\mathcal{S^{*}_{\sum}}}^{\eta}_{q}(\mu,\lambda;\phi)$
if the following subordinations are satisfied:
\begin{equation}
\left
\{\frac{2z[\widetilde{\mathcal{D}}_{q}(\mathcal{J}^{\eta}_{q}f)(z)]^{\lambda}}{\mathcal{J}^{\eta}_{q}f(z)-\mathcal{J}^{\eta}_{q}f(-z)}\right\}^{\mu}\left
\{\frac{2\{\widetilde{\mathcal{D}_{q}}[z\widetilde{\mathcal{D}}_{q}(\mathcal{J}^{\eta}_{q}f)(z)]\}^{\lambda}}{\widetilde{\mathcal{D}}_{q}[\mathcal{J}^{\eta}_{q}f(z)-\mathcal{J}^{\eta}_{q}f(-z)]}\right\}^{1-\mu}\prec\phi(z)
\end{equation}
and
\begin{equation}
\left
\{\frac{2w[\widetilde{\mathcal{D}}_{q}(\mathcal{J}^{\eta}_{q}g)(w)]^{\lambda}}{\mathcal{J}^{\eta}_{q}g(w)-\mathcal{J}^{\eta}_{q}g(-w)}\right\}^{\mu}\left
\{\frac{2\{\widetilde{\mathcal{D}}_{q}[w\widetilde{\mathcal{D}}_{q}(\mathcal{J}^{\eta}_{q}g)(w)]\}^{\lambda}}{\widetilde{\mathcal{D}}_{q}[\mathcal{J}^{\eta}_{q}g(w)-\mathcal{J}^{\eta}_{q}g(-w)]}\right\}^{1-\mu}\prec\phi(w)
\end{equation}
for $z, w\in\Delta$, where $\mu\geq0$ and $\lambda\geq1$, and the
function $g$ is the inverse of $f$ and given by (\ref{equ1.2}).
\end{definition}

\begin{remark}\label{rem1} If we remove the generalized  Bernardi integral operator, then the class
$\widetilde{\mathcal{S^{*}_{\sum}}}^{\eta}_{q}(\mu,\lambda;\phi)$
is exactly
$\widetilde{\mathcal{S^{*}_{\sum}}}_{q}(\mu,\lambda;\phi)$.
Furthermore, the classes
$\widetilde{\mathcal{S^{*}_{\sum}}}^{\eta}_{q}(1,\lambda;\phi)$
and $\mathcal{S^{*}_{\sum}}^{b}_{\eta}(0,\lambda;\phi)$ reduce to
$\widetilde{\mathcal{S^{*}_{\sum}}}_{q}(\lambda;\phi)$ and
$\widetilde{\mathcal{C_{\sum}}}_{q}(\lambda;\phi)$, respectively;
refer to \cite{ES18} for $\mathcal{S^{*}_{\sum}}(\lambda;\phi)$
and $\mathcal{C_{\sum}}(\lambda;\phi)$ departed from the symmetric
$q$-derivative operator,.
\end{remark}

\begin{lemma} [\cite{Du83,Go83}]\label{lem1} Let $\mathcal{P}$ be the class of all analytic functions
$h(z)$  being of the form
\begin{equation*}
h(z)=1+\sum^{\infty}_{n=1}p_{n}z^{n},~~(z\in\Delta)
\end{equation*}
satisfying $\Re h(z)>0$ and $h(0)=1$. Then the sharp estimates
$\mid p_{n}\mid\leq2(n\in\mathbb{N})$ hold. Specially, the
equality is true for all $n$ when
\begin{equation*}
h(z)=\frac{1+z}{1-z}=1+\sum^{\infty}_{n=1}2z^{n}.
\end{equation*}
\end{lemma}

\begin{lemma}[\cite{GS58}]\label{lem2} If $h(z)\in\mathcal{P}$, then the equalities
\begin{equation*}
2p_{2}=p^{2}_{1}+x(4-p^{2}_{1})
\end{equation*}
and
\begin{equation*}
4p_{3}=p^{3}_{1}+2p_{1}(4-p^{2}_{1})x-p_{1}(4-p^{2}_{1})x^{2}+2(4-p^{2}_{1})(1-\vert
x\vert^{2})s
\end{equation*}
hold for some $x,s$ with $\vert x\vert\leq1$ and $\vert
s\vert\leq1$.
\end{lemma}

\section{Functional estimates for the function class
$\widetilde{\mathcal{S^{*}_{\sum}}}_{q}(\lambda;\phi)$}

\hskip\parindent Define the functions $s$ and $t$ in $\mathcal{P}$
by
\begin{equation}\label{equ2.1}
s(z)=\frac{1+u(z)}{1-u(z)}=1+\sum^{\infty}_{n=1}c_{n}z^{n}\quad\text{and}\quad
t(w)=\frac{1+v(w)}{1-v(w)}=
1+\sum^{\infty}_{n=1}d_{n}w^{n},~~(z,w\in\Delta).
\end{equation}
Naturally, from (\ref{equ2.1}) we yield that
\begin{equation}\label{equ2.2}
u(z)=\frac{s(z)-1}{s(z)+1}=\frac{c_{1}}{2}z+\frac{1}{2}\left(c_{2}-\frac{c^{2}_{1}}{2}\right)z^{2}+\frac{1}{2}\left(c_{3}-c_{1}c_{2}+\frac{c^{3}_{1}}{4}\right)z^{3}+\ldots,
~~(z\in\Delta)
\end{equation}
and
\begin{equation}\label{equ2.3}
v(w)=\frac{t(w)-1}{t(w)+1}=
\frac{d_{1}}{2}w+\frac{1}{2}\left(d_{2}-\frac{d^{2}_{1}}{2}\right)w^{2}+\frac{1}{2}\left(d_{3}-d_{1}d_{2}+\frac{d^{3}_{1}}{4}\right)w^{3}\ldots,~~(w\in\Delta).
\end{equation}
On the basis of (\ref{equ2.4}) and (\ref{equ2.2}-\ref{equ2.3}) it infers
that
\begin{eqnarray}\label{equ2.5}
\phi(u(z))&=&1+\frac{1}{2}E_{1}c_{1}z+\left[\frac{1}{2}E_{1}\left(c_{2}-\frac{c^{2}_{1}}{2}\right)+\frac{1}{4}E_{2}c^{2}_{1}\right]z^{2}\notag\\
&+&\left[\frac{1}{2}E_{1}\left(c_{3}-c_{1}c_{2}+\frac{c^{3}_{1}}{4}\right)+\frac{1}{4}E_{2}c_{1}(c_{2}-\frac{c^{2}_{1}}{2})+\frac{1}{8}E_{3}c^{3}_{1}\right]z^{3}+\ldots,~~(z\in\Delta)
\end{eqnarray}
and
\begin{eqnarray}\label{equ2.6}
\phi(v(w))&=&1+\frac{1}{2}E_{1}d_{1}w+\left[\frac{1}{2}E_{1}\left(d_{2}-\frac{d^{2}_{1}}{2}\right)+\frac{1}{4}E_{2}d^{2}_{1}\right]w^{2}\notag\\
&+&\left[\frac{1}{2}E_{1}\left(d_{3}-d_{1}d_{2}+\frac{d^{3}_{1}}{4}\right)+\frac{1}{4}E_{2}d_{1}(d_{2}-\frac{d^{2}_{1}}{2})+\frac{1}{8}E_{3}d^{3}_{1}\right]w^{3}+\ldots,~~(w\in\Delta).
\end{eqnarray}

In this section we are only devoted to study the second Hankel determinant for
the class $\widetilde{\mathcal{S^{*}_{\sum}}}_{q}(\lambda;\phi)$
and write the theorem below.

\begin{theorem}\label{thm1}
If $f(z)$ given by (\ref{equ1.1}) belongs to the class
$\widetilde{\mathcal{S^{*}_{\sum}}}_{q}(\lambda;\phi)$, then
$$
\vert a_{2}a_{4}-a^{2}_{3}\vert\leq\left\{\begin{array}{ll}
\frac{E^{2}_{1}}{(3\lambda-1)^{2}},~~&\mbox{if}~~Q\leq0,~P\leq-\frac{Q}{4},\\
16P+4Q+R,~~ &\mbox{if}~~Q\geq0,~P\geq-\frac{Q}{8}~\text{or}~Q\leq0,~P\geq-\frac{Q}{4},\\
\frac{4PQ-Q^{2}}{4P},~~&\mbox{if}~~Q\geq0,~P\geq-\frac{Q}{8}.
\end{array}\right.
$$
where
\begin{eqnarray*}
P&=&[6(\vert
E_{1}-E_{2}+E_{3}\vert-2E_{1})\lambda^{3}\widetilde{[2]}^{3}_{q}+E^{3}_{1}\{\lambda(\lambda-1)
\vert\lambda-2\vert\widetilde{[2]}^{3}_{q}\notag\\
&+&3\lambda(\lambda-1)\widetilde{[2]}_{q}\widetilde{[3]}_{q}
+6\lambda(\widetilde{[4]}_{q}-\widetilde{[2]}_{q})+15(\lambda-1)\}
]\frac{E^{}_{1}}{96\lambda^{5}\widetilde{[2]}^{4}_{q}\widetilde{[4]}_{q}},
\end{eqnarray*}
\begin{equation*}
Q=\frac{E_{1}(5E_{1}+2\vert
2E_{1}-E_{2}\vert)}{8\lambda^{2}\widetilde{[2]}_{q}\widetilde{[4]}_{q}}+\frac{5
E^{2}_{1}}{4\lambda^{3}\widetilde{[4]}_{q}\widetilde{[2]}^{2}_{q}(\lambda\widetilde{[3]}_{q}-1)}+\frac{E^{3}_{1}}{2\lambda^{2}(3\widetilde{[3]}_{q}-1)\widetilde{[2]}_{q}}
\end{equation*}
and
\begin{equation*}
R=\frac{E^{2}_{1}}{(\lambda\widetilde{[3]}_q-1)^{2}}.
\end{equation*}
\end{theorem}
\begin{proof}[\bf Proof.]  Assume that
$f(z)\in\mathcal{S^{*}_{\sum}}(\lambda;\phi)$. Then, according to
Remark \ref{rem1} and Lemma \ref{lem1} there exist two analytic
functions $u(z)$ and $v(w)\in\mathcal{P}$ so that
\begin{equation}\label{equ2.7}
\frac{2z[\widetilde{\mathcal{D}}_{q}f(z)]^{\lambda}}{f(z)-f(-z)}=\phi(u(z))
\end{equation}
and
\begin{equation}\label{equ2.8}
\frac{2w[\widetilde{\mathcal{D}}_{q}g(w)]^{\lambda}}{g(w)-g(-w)}=\phi(v(w)).
\end{equation}
Clearly, the left hand sides of (\ref{equ2.7}) and (\ref{equ2.8})
can be expanded into the following forms:
\begin{eqnarray}\label{equ2.9}
\frac{2z[\widetilde{\mathcal{D}}_{q}f(z)]^{\lambda}}{f(z)-f(-z)}&=&1+\lambda\widetilde{[2]}_{q}a_{2}z+\left[(\lambda\widetilde{[3]}_{q}-1)a_{3}+\frac{\lambda(\lambda-1)}{2}\widetilde{[2]}^{2}_{q}a^{2}_{2}\right]z^{2}\notag\\
 &+&\left\{\lambda \widetilde{[4]}_{q}a_{4}+\lambda\left(\frac{\lambda-1}{2}\widetilde{[3]}_{q}-1\right)\widetilde{[2]}_{q}a_{2}a_{3}+\frac{\lambda(\lambda-1)(\lambda-2)}{6}\widetilde{[2]}^{3}_{q}a^{3}_{2}\right\}z^{3}+\ldots
\end{eqnarray}
and
\begin{eqnarray}\label{equ2.10}
\frac{2w[\widetilde{\mathcal{D}}_{q}g(w)]^{\lambda}}{g(w)-g(-w)}&=&1-\lambda\widetilde{[2]}_{q}
a_{2}w+\left[(\lambda\widetilde{[3]}_{q}-1)(2a^{2}_{2}-a_{3})
+\frac{\lambda(\lambda-1)}{2}\widetilde{[2]}^{2}_{q}a^{2}_{2}\right]w^{2}\notag\\
&-&[\lambda(5a^{3}_{2}-5a_{2}a_{3}+a_{4})\widetilde{[4]}_{q}+\lambda\left(\frac{\lambda-1}{2}\widetilde{[3]}_{q}-1\right)\widetilde{[2]}_{q}a_{2}(2a^{2}_{2}-a_{3})\notag\\
&+&\frac{\lambda(\lambda-1)(\lambda-2)}{6}\widetilde{[2]}_{q}^{3}a^{3}_{2}]w^{3}+\ldots
\end{eqnarray}
Hence, from (\ref{equ2.5}-\ref{equ2.6}) and
(\ref{equ2.7}-\ref{equ2.10}) we get that
\begin{equation}\label{equ2.11}
\lambda\widetilde{[2]}_{q}a_{2}=\frac{1}{2}E_{1}c_{1},
\end{equation}
\begin{equation}\label{equ2.12}
(\lambda\widetilde{[3]}_{q}-1)a_{3}+\frac{\lambda(\lambda-1)}{2}\widetilde{[2]}^{2}_{q}a^{2}_{2}=\frac{1}{2}E_{1}\left(c_{2}-\frac{c^{2}_{1}}{2}\right)+\frac{1}{4}E_{2}c^{2}_{1},
\end{equation}
\begin{eqnarray}\label{equ2.13}
&&\lambda \widetilde{[4]}_{q}a_{4}+\lambda\left(\frac{\lambda-1}{2}\widetilde{[3]}_{q}-1\right)\widetilde{[2]}_{q}a_{2}a_{3}+\frac{\lambda(\lambda-1)(\lambda-2)}{6}\widetilde{[2]}^{3}_{q}a^{3}_{2}\notag\\
&=&\frac{1}{2}E_{1}\left(c_{3}-c_{1}c_{2}+\frac{c^{3}_{1}}{4}\right)+\frac{1}{4}E_{2}c_{1}\left(c_{2}-\frac{c^{2}_{1}}{2}\right)+\frac{1}{8}E_{3}c^{3}_{1},
\end{eqnarray}
\begin{equation}\label{equ2.14}
-\lambda\widetilde{[2]}_{q} a_{2}=\frac{1}{2}E_{1}d_{1},
\end{equation}
\begin{equation}\label{equ2.15}
(\lambda\widetilde{[3]}_{q}-1)(2a^{2}_{2}-a_{3})
+\frac{\lambda(\lambda-1)}{2}\widetilde{[2]}^{2}_{q}a^{2}_{2}=\frac{1}{2}E_{1}\left(d_{2}-\frac{d^{2}_{1}}{2}\right)+\frac{1}{4}E_{2}d^{2}_{1}
\end{equation}
and
\begin{eqnarray}\label{equ2.16}
-\lambda(5a^{3}_{2}-5a_{2}a_{3}+a_{4})\widetilde{[4]}_{q}&-&\lambda\left(\frac{\lambda-1}{2}\widetilde{[3]}_{q}-1\right)\widetilde{[2]}_{q}a_{2}(2a^{2}_{2}-a_{3})
-\frac{\lambda(\lambda-1)(\lambda-2)}{6}\widetilde{[2]}_{q}^{3}a^{3}_{2}\notag\\
&=&\frac{1}{2}E_{1}\left(d_{3}-d_{1}d_{2}+\frac{d^{3}_{1}}{4}\right)+\frac{1}{4}E_{2}d_{1}\left(d_{2}-\frac{d^{2}_{1}}{2}\right)+\frac{1}{8}E_{3}d^{3}_{1}.
\end{eqnarray}
Further, by (\ref{equ2.11}) and (\ref{equ2.12}) we derive that
\begin{equation}\label{equ2.17}
a_{2}=\frac{E_{1}c_{1}}{2\lambda\widetilde{[2]}_{q}}=-\frac{E_{1}d_{1}}{2\lambda\widetilde{[2]}_{q}}
\end{equation}
such that
\begin{equation}\label{equ2.18}
c_{1}=-d_{1}
\end{equation}
and
\begin{equation}\label{equ2.19}
E^{2}_{1}(c^{2}_{1}+d^{2}_{1})=8\lambda^{2}\widetilde{[2]}^{2}_{q}a^{2}_{2}.
\end{equation}
Making use of (\ref{equ2.12}) and (\ref{equ2.15}) we also obtain that
\begin{equation}\label{equ2.20}
2(\lambda\widetilde{[3]}_{q}-1)(a_{3}-a^{2}_{2})=\frac{1}{2}E_{1}(c_{2}-d_{2}).
\end{equation}
Therefore
\begin{equation}\label{equ2.21}
a_{3}=\frac{E_{1}(c_{2}-d_{2})}{4(\lambda\widetilde{[3]}_{q}-1)}+\frac{E^{2}_{1}c^{2}_{1}}{4\lambda^{2}\widetilde{[2]}^{2}_{q}}.
\end{equation}
In addition, by (\ref{equ2.13}) and (\ref{equ2.16}) we conclude that

\begin{eqnarray}\label{equ2.22}
a_{4}&=&-\frac{5}{2\widetilde{[4]}_{q}}a_{2}(a^{2}_{2}-a_{3})-\frac{\widetilde{[2]}_{q}}{6\widetilde{[4]}_{q}}\{(\lambda-1)[3\widetilde{[3]}_{q}+(\lambda-2)\widetilde{[2]}^{2}_{q}]-6\}a^{3}_{2}\notag\\
&+&\frac{E_{1}}{4\lambda\widetilde{[4]}_{q}}(c_{3}-d_{3})+\frac{c_{1}}{8\lambda\widetilde{[4]}_{q}}(E_{2}-2E_{1})(c_{2}+d_{2})+\frac{c^{3}_{1}}{8\lambda\widetilde{[4]}_{q}}(E_{1}-E_{2}+E_{3}),
\end{eqnarray}
then we show from (\ref{equ2.17}) and
(\ref{equ2.21}-\ref{equ2.22}) that
\begin{eqnarray}\label{equ2.23}
a_{2}a_{4}-a^{2}_{3}&=&-\frac{5}{2\widetilde{[4]}_{q}}a^{2}_{2}(a^{2}_{2}-a_{3})-\frac{\widetilde{[2]}_{q}}{6\widetilde{[4]}_{q}}\{(\lambda-1)[3\widetilde{[3]}_{q}+(\lambda-2)\widetilde{[2]}^{2}_{q}]-6\}a^{4}_{2}\notag\\
&+&\frac{E_{1}a_{2}}{4\lambda\widetilde{[4]}_{q}}(c_{3}-d_{3})+\frac{c_{1}a_{2}}{8\lambda\widetilde{[4]}_{q}}(E_{2}-2E_{1})(c_{2}+d_{2})+\frac{c^{3}_{1}a_{2}}{8\lambda\widetilde{[4]}_{q}}(E_{1}-E_{2}+E_{3})-a^{2}_{3}\notag\\
&=&\frac{E_{1}c^{4}_{1}}{96\lambda^{4}\widetilde{[2]}^{4}_{q}\widetilde{[4]}_{q}}[-E^{3}_{1}\{(\lambda-1)
(\lambda-2)\widetilde{[2]}^{3}_{q}+3(\lambda-1)\widetilde{[2]}_{q}\widetilde{[3]}_{q}-6(\widetilde{[2]}_{q}-\widetilde{[4]}_{q})\}\notag\\
&+&6(E_{1}-E_{2}+E_{3})\lambda^{2}\widetilde{[2]}^{3}_{q}]+\frac{E^{2}_{1}c_{1}(c_{3}-d_{3})}{8\lambda^{2}\widetilde{[2]}_{q}\widetilde{[4]}_{q}}
+\frac{E_{1}(E_{2}-2E_{1})c^{2}_{1}(c_{2}+d_{2})}{16\lambda^{2}\widetilde{[2]}_{q}\widetilde{[4]}_{q}}\notag\\
&-&\frac{E^{2}_{1}(c_{2}-d_{2})^{2}}{16(\lambda\widetilde{[3]}_{q}-1)^{2}}
+\frac{(5-4\widetilde{[4]}_{q})E^{3}_{1}c^{2}_{1}(c_{2}-d_{2})}{32\lambda^{2}\widetilde{[2]}^{2}_{q}\widetilde{[4]}_{q}(\lambda\widetilde{[3]}_{q}-1)}.
\end{eqnarray}
From Lemma \ref{lem2} and (\ref{equ2.18}), it derives that
\begin{equation}\label{equ2.24}
c_{2}-d_{2}=\frac{1}{2}(4-c^{2}_{1})(x-y) \quad\text{and}\quad
c_{2}+d_{2}=c^{2}_{1}+\frac{1}{2}(4-c^{2}_{1})(x+y),
\end{equation}
and
\begin{equation}\label{equ2.25}
c_{3}-d_{3}=\frac{c^{3}_{1}}{2}+\frac{c_{1}}{2}(4-c^{2}_{1})(x+y)-\frac{c_{1}}{4}(4-c^{2}_{1})(x^{2}+y^{2})+\frac{1}{2}(4-c^{2}_{1})[(1-\vert
x\vert^{2})s-(1-\vert y\vert^{2})t]
\end{equation}
for some $x,y,s,t\in [-1,1]$. Together with (\ref{equ2.24}) and
(\ref{equ2.25}), (\ref{equ2.23}) follows that
\begin{eqnarray}\label{equ2.26}
a_{2}a_{4}-a^{2}_{3}&=&\frac{-E^{}_{1}c^{4}_{1}}{96\lambda^{4}\widetilde{[2]}^{4}_{q}\widetilde{[4]}_{q}}[E^{3}_{1}\{(\lambda-1)
(\lambda-2)\widetilde{[2]}^{3}_{q}+3(\lambda-1)\widetilde{[2]}_{q}\widetilde{[3]}_{q}
-6(\widetilde{[2]}_{q}-\widetilde{[4]}_{q})\}\notag\\
&-&6E_{3}\lambda^{2}\widetilde{[2]}^{3}_{q}]+\frac{E_{1}E_{2}c^{2}_{1}(4-c^{2}_{1})(x+y)}{32\lambda^{2}\widetilde{[2]}_{q}\widetilde{[4]}_{q}}\notag\\
&-&\frac{E^{2}_{1}c^{2}_{1}(4-c^{2}_{1})}{32\lambda^{2}\widetilde{[2]}_{q}\widetilde{[4]}_{q}}(
x^{2}+
y^{2})+\frac{E^{2}_{1}c_{1}(4-c^{2}_{1})}{16\lambda^{2}\widetilde{[2]}_{q}\widetilde{[4]}_{q}}[(1-
x^{2})s-(1-y^{2})t]\notag\\
&-&\frac{E^{2}_{1}(4-c^{2}_{1})^{2}(x-y)^{2}}{64(\lambda\widetilde{[3]}_{q}-1)^{2}}
+\frac{(5-4\widetilde{[4]}_{q})E^{3}_{1}c^{2}_{1}(4-c^{2}_{1})(x-y)}{64\lambda^{2}\widetilde{[2]}^{2}_{q}\widetilde{[4]}_{q}(\widetilde{\lambda[3]}_{q}-1)}.
\end{eqnarray}
Then
\begin{eqnarray*}
\vert
a_{2}a_{4}-a^{2}_{3}\vert&\leq&\frac{E^{}_{1}c^{4}_{1}}{96\lambda^{4}\widetilde{[2]}^{4}_{q}\widetilde{[4]}_{q}}[E^{3}_{1}\{(\lambda-1)
\vert\lambda-2\vert\widetilde{[2]}^{3}_{q}+3(\lambda-1)\widetilde{[2]}_{q}\widetilde{[3]}_{q}
+6(\widetilde{[4]}_{q}-\widetilde{[2]}_{q})\}\notag\\
&+&6\vert E_{3}\vert\lambda^{2}
\widetilde{[2]}^{3}_{q}]+\frac{E_{1}\vert E_{2}\vert
c^{2}_{1}(4-c^{2}_{1})(\vert x\vert+\vert
y\vert)}{32\lambda^{2}\widetilde{[2]}_{q}\widetilde{[4]}_{q}}\notag\\
&+&\frac{E^{2}_{1} c^{2}_{1}(4-c^{2}_{1})(\vert x\vert^{2}+ \vert
y\vert^{2})}{32\lambda^{2}\widetilde{[2]}_{q}\widetilde{[4]}_{q}}
+\frac{E^{2}_{1}c_{1}(4-c^{2}_{1})}{16\lambda^{2}\widetilde{[2]}_{q}\widetilde{[4]}_{q}}
[(1-\vert x\vert^{2})+(1-\vert y\vert^{2})]\notag\\
&+&\frac{E^{2}_{1}(4-c^{2}_{1})^{2}(\vert x\vert+\vert
y\vert)^{2}}{64(\lambda\widetilde{[3]}_{q}-1)^{2}}
+\frac{\vert5-4\widetilde{[4]}_{q}\vert
E^{3}_{1}c^{2}_{1}(4-c^{2}_{1})(\vert x\vert+\vert
y\vert)}{64\lambda^{2}\widetilde{[2]}^{2}_{q}\widetilde{[4]}_{q}(\lambda\widetilde{[3]}_{q}-1)}\notag\\
&=&\frac{E^{}_{1}c^{4}_{1}}{96\lambda^{4}\widetilde{[2]}^{4}_{q}\widetilde{[4]}_{q}}[E^{3}_{1}\{(\lambda-1)
\vert\lambda-2\vert\widetilde{[2]}^{3}_{q}+3(\lambda-1)\widetilde{[2]}_{q}\widetilde{[3]}_{q}
+6(\widetilde{[4]}_{q}-\widetilde{[2]}_{q})\}\notag\\&+&6\vert E_{3}\vert\lambda^{2}\widetilde{[2]}^{3}_{q}]+\frac{E^{2}_{1}c_{1}(4-c^{2}_{1})}{8\lambda^{2}\widetilde{[2]}_{q}\widetilde{[4]}_{q}}\notag\\
&+&\left\{\frac{E_{1}\vert E_{2}\vert
}{32\lambda^{2}\widetilde{[2]}_{q}\widetilde{[4]}_{q}}+\frac{\vert5-4\widetilde{[4]}_{q}\vert
E^{3}_{1}}{64\lambda^{2}\widetilde{[2]}^{2}_{q}\widetilde{[4]}_{q}(\lambda\widetilde{[3]}_{q}-1)}\right\}
c^{2}_{1}(4-c^{2}_{1})(\vert x\vert+\vert y\vert)\notag\\
&-&\frac{E^{2}_{1}c_{1}(2-c_{1})(4-c^{2}_{1})}{32\lambda^{2}\widetilde{[2]}_{q}\widetilde{[4]}_{q}}(\vert
x\vert^{2}+\vert
y\vert^{2})+\frac{E^{2}_{1}(4-c^{2}_{1})^{2}(\vert x\vert+\vert
y\vert)^{2}}{64(\lambda\widetilde{[3]}_{q}-1)^{2}}.
\end{eqnarray*}
Since $\vert c_{1}\vert\leq2$, we assume that $c=c_{1}\in [0,2]$.
Letting $\xi=\vert x\vert$ and $\zeta=\vert y\vert$, we remark
that
\begin{equation}
\vert
a_{2}a_{4}-a^{2}_{3}\vert\leq\mathcal{F}_{1}+(\xi+\zeta)\mathcal{F}_{2}+(\xi^{2}+\zeta^{2})\mathcal{F}_{3}+(\xi+\zeta)^{2}\mathcal{F}_{4}=:\Pi(\xi,\zeta)
\end{equation}
in the closed square $\mathcal{D}=\{(\xi,\zeta):
0\leq\xi\leq1,0\leq\zeta\leq1)\}$, where
\begin{eqnarray*}
\mathcal{F}_{1}=\mathcal{F}_{1}(c)&:=&\frac{E^{}_{1}c^{4}}{96\lambda^{4}\widetilde{[2]}^{4}_{q}\widetilde{[4]}_{q}}[E^{3}_{1}\{(\lambda-1)
\vert\lambda-2\vert\widetilde{[2]}^{3}_{q}+3(\lambda-1)\widetilde{[2]}_{q}\widetilde{[3]}_{q}
+6(\widetilde{[4]}_{q}-\widetilde{[2]}_{q})\}\notag\\&+&6\vert
E_{3}\vert\lambda^{2}\widetilde{[2]}^{3}_{q}]+\frac{E^{2}_{1}c(4-c^{2})}{8\lambda^{2}\widetilde{[2]}_{q}\widetilde{[4]}_{q}}\geq0,
\end{eqnarray*}
\begin{equation*}
\mathcal{F}_{2}=\mathcal{F}_{2}(c):=\left\{\frac{E_{1}\vert
E_{2}\vert}{32\lambda^{2}\widetilde{[2]}_{q}\widetilde{[4]}_{q}}+\frac{\vert5-4\widetilde{[4]}_{q}\vert
E^{3}_{1}}{64\lambda^{2}\widetilde{[2]}^{2}_{q}\widetilde{[4]}_{q}(\lambda\widetilde{[3]}_{q}-1)}\right\}c^{2}(4-c^{2})\geq0,
\end{equation*}
\begin{equation*}
\mathcal{F}_{3}=\mathcal{F}_{3}(c):=-\frac{E^{2}_{1}c(2-c)(4-c^{2})}{32\lambda^{2}\widetilde{[2]}_{q}\widetilde{[4]}_{q}}\leq0\quad
and\quad\mathcal{F}_{4}=\mathcal{F}_{4}(c):=\frac{E^{2}_{1}(4-c^{2})^{2}}{64(\lambda\widetilde{[3]}_{q}-1)^{2}}\geq0.
\end{equation*}

For $c\in[0,2]$, we will maximize $\Pi(\xi,\zeta)$ by the sign of
$\Pi_{\xi\xi}\Pi_{\zeta\zeta}-(\Pi_{\xi\zeta})^{2}$ in the closed
square $\mathcal{D}=\{(\xi,\zeta):
0\leq\xi\leq1,0\leq\zeta\leq1)\}$.

Consider $c\in(0,2)$. By the simple calculation, we observe that
$\mathcal{F}_{3}<0$ and $\mathcal{F}_{3}+2\mathcal{F}_{4}>0$ for
$c\in(0,2)$. Then, we know that
$\Pi_{\xi\xi}\Pi_{\zeta\zeta}-(\Pi_{\xi\zeta})^{2}<0$ such that
$\Pi$ doesn't take the local maximum in the interior the square
$\mathcal{D}$.

For $\xi=0$ and $0\leq\zeta\leq1$ (or, for $\zeta=0$ and
$0\leq\xi\leq1$), we see that
\begin{equation*}
\Pi(0,\zeta)=H(\zeta)=\mathcal{F}_{1}+\mathcal{F}_{2}\zeta+(\mathcal{F}_{3}+\mathcal{F}_{4})\zeta^{2}.
\end{equation*}

Case (I) If $\mathcal{F}_{3}+\mathcal{F}_{4}\geq0$, then
$H'(\zeta)=\mathcal{F}_{2}+2(\mathcal{F}_{3}+\mathcal{F}_{4})\zeta
>0$ for $0<\zeta<1$ and $c\in(0,2)$, and $H(\zeta)$ is an
increasing function in $[0,1]$. Hence, $H(\zeta)$ reaches the
maximum via
\begin{equation*}
\max
H(\zeta)=H(1)=\mathcal{F}_{1}+\mathcal{F}_{2}+\mathcal{F}_{3}+\mathcal{F}_{4}.
\end{equation*}

Case (II) If $\mathcal{F}_{3}+\mathcal{F}_{4}<0$, we study two
versions for the critical point
$\widetilde{\zeta}=\frac{-\mathcal{F}_{2}}{2(\mathcal{F}_{3}+\mathcal{F}_{4})}$.
When $\widetilde{\zeta}>1$, we infer that
$0\leq-(\mathcal{F}_{3}+\mathcal{F}_{4})\leq\mathcal{F}_{2}+\mathcal{F}_{3}+\mathcal{F}_{4}$
so that
$H(0)\leq\mathcal{F}_{1}+\mathcal{F}_{2}+\mathcal{F}_{3}+\mathcal{F}_{4}=
H(1)$. Inversely, when $\widetilde{\zeta}\leq1$, we have that
\begin{equation*}
H(1)=\mathcal{F}_{1}+\frac{1}{2}\mathcal{F}_{2}+\frac{1}{2}(\mathcal{F}_{2}+2\mathcal{F}_{3}+2\mathcal{F}_{4})\leq\mathcal{F}_{1}+\frac{1}{2}\mathcal{F}_{2}
\end{equation*}
and
\begin{equation*}
H(0)\leq
H(\widetilde{\zeta})=\mathcal{F}_{1}+\frac{-\mathcal{F}^{2}_{2}}{4(\mathcal{F}_{3}+\mathcal{F}_{4})}\leq\mathcal{F}_{1}+\frac{1}{2}\mathcal{F}_{2}.
\end{equation*}
As for $c=0$ or $c=2$, we see that
\begin{equation*}
\Pi(\xi,\zeta)=(\xi+\zeta)^{2}\mathcal{F}_{4},\mathcal{F}_{1}=\mathcal{F}_{2}=\mathcal{F}_{3}=0~\text{and}~\mathcal{T}_{4}>0
\end{equation*}
or
\begin{equation*}
\Pi(\xi,\zeta)=\mathcal{F}_{1}+(\xi+\zeta)^{2}\mathcal{F}_{4},\mathcal{F}_{2}=\mathcal{F}_{3}=0~
\text{and}~\mathcal{F}_{1},\mathcal{F}_{4}>0.
\end{equation*}
When $\xi=\zeta=1$, we note that $F$ reaches the corresponding
maximum $4\mathcal{T}_{4}$ or $\mathcal{F}_{1}+4\mathcal{F}_{4}$.

For $\xi=1$ and $0\leq\zeta\leq1$ (or, for $\zeta=1$ and
$0\leq\xi\leq1$), we derive that
\begin{equation*}
\Pi(1,\zeta)=:G(\zeta)=\mathcal{F}_{1}+\mathcal{F}_{2}+\mathcal{F}_{3}+\mathcal{F}_{4}+(\mathcal{F}_{2}+2\mathcal{F}_{4})\zeta+(\mathcal{F}_{3}+\mathcal{F}_{4})\zeta^{2}.
\end{equation*}
Case (III) If $\mathcal{F}_{3}+\mathcal{F}_{4}\geq0$, then
$G'(\zeta)=\mathcal{F}_{2}+2\mathcal{F}_{4}+2(\mathcal{F}_{3}+\mathcal{F}_{4})\zeta
>0$ for $0<\zeta<1$ and $c\in(0,2)$, and so $G(\zeta)$ is also an
increasing function in $[0,1]$. Therefore, the maximum of
$G(\zeta)$ is
\begin{equation*}
\max
G(\zeta)=G(1)=\mathcal{F}_{1}+2\mathcal{F}_{2}+2\mathcal{F}_{3}+4\mathcal{F}_{4}.
\end{equation*}
Case (IV) If $\mathcal{F}_{3}+\mathcal{F}_{4}<0$, then there have
two versions with respect to the critical point
$\widetilde{\zeta}=-\frac{\mathcal{F}_{2}+2\mathcal{F}_{4}}{2(\mathcal{F}_{3}+\mathcal{F}_{4})}$.
When $\widetilde{\zeta}>1$, we know that
$0\leq-(\mathcal{F}_{3}+\mathcal{F}_{4})\leq\mathcal{F}_{2}+\mathcal{F}_{3}+3\mathcal{F}_{4}$
so that
$G(0)\leq\mathcal{F}_{1}+2\mathcal{F}_{2}+2\mathcal{F}_{3}+4\mathcal{F}_{4}=
G(1)$. Oppositely, when $\widetilde{\zeta}\leq1$,
\begin{equation*}
G(0)\leq
G(\widetilde{\zeta})=\mathcal{F}_{1}+\mathcal{F}_{2}+\mathcal{F}_{3}+\mathcal{F}_{4}-\frac{(\mathcal{F}_{2}+2\mathcal{F}_{4})^{2}}{4(\mathcal{F}_{3}+\mathcal{F}_{4})}
\leq\mathcal{F}_{1}+\frac{3}{2}\mathcal{F}_{2}+\mathcal{F}_{3}+2\mathcal{F}_{4}\leq
G(1).
\end{equation*}
All in all, by all the above cases, the maximum of $\Pi$ in the
closed square $\mathcal{D}=\{(\xi,\zeta):
0\leq\xi\leq1,0\leq\zeta\leq1)\}$ occurs at $\xi=1$ and $\zeta=1$,
that is to say,
\begin{equation*}
\max_{(\xi,\zeta)\in\mathcal{D}}\Pi(\xi,\zeta)=\Pi(1,1).
\end{equation*}

Consider the function $\mathcal{G}:[0,2]\rightarrow\mathbb{R}$ by
\begin{equation*}
\mathcal{G}(c)=\max_{(\xi,\zeta)\in\mathcal{D}}\Pi(\xi,\zeta)=\Pi(1,1)=\mathcal{F}_{1}+2\mathcal{F}_{2}+2\mathcal{F}_{3}+4\mathcal{F}_{4}.
\end{equation*}
Then
\begin{eqnarray*}
\mathcal{G}(c)&=&\frac{E_{1}}{96\lambda^{4}\widetilde{[2]}^{4}_{q}\widetilde{[4]}_{q}}[E^{3}_{1}\{(\lambda-1)
\vert\lambda-2\vert\widetilde{[2]}^{3}_{q}+3(\lambda-1)\widetilde{[2]}_{q}\widetilde{[3]}_{q}
+6(\widetilde{[4]}_{q}-\widetilde{[2]}_{q})\}\notag\\
&+&6(\vert E_{3}\vert-E_{1}-\vert
E_{2}\vert)\lambda^{2}\widetilde{[2]}^{3}_{q}-\frac{3\vert5-4\widetilde{[4]}_{q}\vert
E^{2}_{1}\lambda^{2}\widetilde{[2]}^{2}_{q}}{\lambda\widetilde{[3]}_{q}-1}+\frac{6E_{1}\lambda^{4}\widetilde{[2]}^{4}_{q}\widetilde{[4]}_{q}}{(\lambda\widetilde{[3]}_{q}-1)^{2}}]c^{4}\notag\\
&+&\left[\frac{E_{1}(E_{1}+\vert
E_{2}\vert)}{4\lambda^{2}\widetilde{[2]}_{q}\widetilde{[4]}_{q}}+\frac{\vert5-4\widetilde{[4]}_{q}\vert
E^{3}_{1}}{8\lambda^{2}\widetilde{[2]}^{2}_{q}\widetilde{[4]}_{q}(\lambda\widetilde{[3]}_{q}-1)}-\frac{
E^{2}_{1}}{2(\lambda\widetilde{[3]}_{q}-1)^{2}}\right]c^{2}\notag\\
&+&\frac{E^{2}_{1}}{(\lambda\widetilde{[3]}_{q}-1)^{2}}.
\end{eqnarray*}
Let $t=c^{2}$. Denote
\begin{eqnarray}\label{equ2.27}
P&=&\frac{E_{1}}{96\lambda^{4}\widetilde{[2]}^{4}_{q}\widetilde{[4]}_{q}}[E^{3}_{1}\{(\lambda-1)
\vert\lambda-2\vert\widetilde{[2]}^{3}_{q}+3(\lambda-1)\widetilde{[2]}_{q}\widetilde{[3]}_{q}
+6(\widetilde{[4]}_{q}-\widetilde{[2]}_{q})\}\notag\\
&+&6(\vert E_{3}\vert-E_{1}-\vert
E_{2}\vert)\lambda^{2}\widetilde{[2]}^{3}_{q}-\frac{3\vert5-4\widetilde{[4]}_{q}\vert
E^{2}_{1}\lambda^{2}\widetilde{[2]}^{2}_{q}}{\lambda\widetilde{[3]}_{q}-1}+\frac{6E_{1}\lambda^{4}\widetilde{[2]}^{4}_{q}\widetilde{[4]}_{q}}{(\lambda\widetilde{[3]}_{q}-1)^{2}}],
\end{eqnarray}
\begin{equation}\label{equ2.28}
Q=\frac{E_{1}(E_{1}+\vert
E_{2}\vert)}{4\lambda^{2}\widetilde{[2]}_{q}\widetilde{[4]}_{q}}+\frac{\vert5-4\widetilde{[4]}_{q}\vert
E^{3}_{1}}{8\lambda^{2}\widetilde{[2]}^{2}_{q}\widetilde{[4]}_{q}(\lambda\widetilde{[3]}_{q}-1)}-\frac{
E^{2}_{1}}{2(\lambda\widetilde{[3]}_{q}-1)^{2}}
\end{equation}
and
\begin{equation}\label{equ2.29}
R=\frac{E^{2}_{1}}{(\lambda\widetilde{[3]}_q-1)^{2}}.
\end{equation}
By folowing the standard computations of the optimal value of
quadratic as follows:
$$
\max_{0\leq t\leq4}(Pt^{2}+Qt+R)=\left\{\begin{array}{ll}
R,~~&\mbox{if}~~Q\leq0,~P\leq-\frac{Q}{4},\\
16P+4Q+R,~~ &\mbox{if}~~Q\geq0,~P\geq-\frac{Q}{8}~\text{or}~Q\leq0,~P\geq-\frac{Q}{4},\\
\frac{4PQ-Q^{2}}{4P},~~&\mbox{if}~~Q\geq0,~P\geq-\frac{Q}{8},
\end{array}\right.
$$
we imply that
$$
\vert a_{2}a_{4}-a^{2}_{3}\vert\leq\left\{\begin{array}{ll}
R,~~&\mbox{if}~~Q\leq0,~P\leq-\frac{Q}{4},\\
16P+4Q+R,~~ &\mbox{if}~~Q\geq0,~P\geq-\frac{Q}{8}~\text{or}~Q\leq0,~P\geq-\frac{Q}{4},\\
\frac{4PQ-Q^{2}}{4P},~~&\mbox{if}~~Q\geq0,~P\geq-\frac{Q}{8},
\end{array}\right.
$$
where $P,Q$ and $R$ are given by (\ref{equ2.27}-\ref{equ2.29}).
\end{proof}

\begin{remark} By the same method as in Theorem \ref{thm1},
further we can study the second Hankel determinants for the
classes $\widetilde{\mathcal{C_{\sum}}}_{q}(\lambda;\phi)$ with
respect to symmetric points, and leave them to the interested
readers. In addition, when $q\rightarrow1_{-}$, there are some
versions of the second Hankel determinants for some function
classes under some appropriate assumption conditions:

(I) For $\phi(z)=z+\sqrt{1+z^{2}}$, refer to \cite{SL18} for the
classes $\widetilde{\mathcal{S^{*}_{\sum}}}_{1_{-}}(1;\phi)$ and
$\widetilde{\mathcal{C_{\sum}}}_{1_{-}}(1;\phi)$ related to a
shell shaped region;

(II) For $\phi(z)=\frac{1+Az}{1+Bz}$, refer to \cite{GM82} for the
classes $\widetilde{\mathcal{S^{*}_{\sum}}}_{1_{-}}(1;\phi)$ and
$\widetilde{\mathcal{C_{\sum}}}_{1_{-}}(1;\phi)$ with respect to a
symmetric or conjugate point.

\end{remark}

\section{Functional estimates for the function class $\widetilde{\mathcal{S^{*}_{\sum}}}^{\eta}_{q}(\mu,\lambda;\phi)$}

\hskip\parindent In this section we are deeply concerned about the functional estimates
for the class
$\widetilde{\mathcal{S^{*}_{\sum}}}^{\eta}_{q}(\mu,\lambda;\phi)$.
Now we state our theorem of coefficient bounds for $\widetilde{\mathcal{S^{*}_{\sum}}}^{\eta}_{q}(\mu,\lambda;\phi)$ as follows.
\begin{theorem}\label{thm2} If $f(z)$ given by (\ref{equ1.1}) is in the class
$\widetilde{\mathcal{S^{*}_{\sum}}}^{\eta}_{q}(\mu,\lambda;\phi)$,
then
\begin{equation}\label{equ3.1}
\mid
a_{2}\mid\leq\min\left\{\frac{E_{1}}{\lambda\vert[\mu-(\mu-1)\widetilde{[2]}_{q}]L_{2}\vert\widetilde{[2]}_{q}},
\sqrt{\frac{2(\vert E_{2}-E_{1}\vert+ E_{1})}{\vert\Omega\vert}},
\frac{E_{1}\sqrt{2E_{1}}}{\sqrt{\vert\Theta\vert}}\right\}
\end{equation}
and
\begin{eqnarray}\label{equ3.2}
\mid a_{3}\mid&\leq&
\frac{E_{1}}{(\lambda\widetilde{[3]}_{q}-1)\vert[\mu-(\mu-1)\widetilde{[3]}_{q}]L_{3}\vert}\notag\\
&+&\min\left\{\frac{E^{2}_{1}}{\lambda^{2}[\mu-(\mu-1)\widetilde{[2]}_{q}]^{2}\widetilde{[2]}^{2}_{q}\vert
L_{2}\vert^{2}}, \frac{\vert
2(E_{2}-E_{1}\vert+E_{1})}{\vert\Omega\vert}\right\},
\end{eqnarray}
where
\begin{equation}\label{equ3.3}
\Omega:=2(\lambda\widetilde{[3]}_{q}-1)[\mu-(\mu-1)\widetilde{[3]}_{q}]L_{3}+\lambda(\lambda[\mu-(\mu-1)\widetilde{[2]}_{q}]^{2}-[\mu-(\mu-1)\widetilde{[2]}^{2}_{q}])\widetilde{[2]}^{2}_{q}L^{2}_{2}
\end{equation}
and
\begin{eqnarray}\label{equ3.4}
\Theta:&=&2(\lambda\widetilde{[3]}_{q}-1)[\mu-(\mu-1)\widetilde{[3]}_{q}]E^{2}_{1}L_{3}\notag\\
&+&\lambda\{\lambda[\mu-(\mu-1)\widetilde{[2]}_{q}]^{2}(E^{2}_{1}+2E_{1}-2E_{2})-[\mu-(\mu-1)\widetilde{[2]}^{2}_{q}]E^{2}_{1}\}\widetilde{[2]}^{2}_{q}L^{2}_{2}.
\end{eqnarray}
\end{theorem}

\begin{proof}[\bf Proof.] Suppose that
$f(z)\in\widetilde{\mathcal{S^{*}_{\sum}}}^{\eta}_{q}(\mu,\lambda;\phi)$.
Then, by Definition \ref{def1} and Lemma \ref{lem1} there exist
two analytic functions $u(z)$ and $v(w)\in\mathcal{P}$ such that
\begin{equation}\label{equ3.5}
\left
\{\frac{2z[\widetilde{\mathcal{D}}_{q}(\mathcal{J}^{\eta}_{q}f)(z)]^{\lambda}}{\mathcal{J}^{\eta}_{q}f(z)-\mathcal{J}^{\eta}_{q}f(-z)}\right\}^{\mu}\left
\{\frac{2\{\widetilde{\mathcal{D}_{q}}[z\widetilde{\mathcal{D}}_{q}(\mathcal{J}^{\eta}_{q}f)(z)]\}^{\lambda}}{\widetilde{\mathcal{D}}_{q}[\mathcal{J}^{\eta}_{q}f(z)-\mathcal{J}^{\eta}_{q}f(-z)]}\right\}^{1-\mu}=\phi(u(z))
\end{equation}
and
\begin{equation}\label{equ3.6}
\left
\{\frac{2w[\widetilde{\mathcal{D}}_{q}(\mathcal{J}^{\eta}_{q}g)(w)]^{\lambda}}{\mathcal{J}^{\eta}_{q}g(w)-\mathcal{J}^{\eta}_{q}g(-w)}\right\}^{\mu}\left
\{\frac{2\{\widetilde{\mathcal{D}}_{q}[w\widetilde{\mathcal{D}}_{q}(\mathcal{J}^{\eta}_{q}g)(w)]\}^{\lambda}}{\widetilde{\mathcal{D}}_{q}[\mathcal{J}^{\eta}_{q}g(w)-\mathcal{J}^{\eta}_{q}g(-w)]}\right\}^{1-\mu}=\phi(v(w))
\end{equation}
By expanding the left half parts of (\ref{equ3.5}) and
(\ref{equ3.6}), it leads to
\begin{eqnarray*}
&&\left
\{\frac{2z[\widetilde{\mathcal{D}}_{q}(\mathcal{J}^{\eta}_{q}f)(z)]^{\lambda}}{\mathcal{J}^{\eta}_{q}f(z)-\mathcal{J}^{\eta}_{q}f(-z)}\right\}^{\mu}=1+\lambda\mu\widetilde{[2]}_{q}L_{2}a_{2}z
\notag\\
 &+&\left\{\mu[(\lambda\widetilde{[3]}_{q}-1)L_{3}a_{3}+\frac{\lambda(\lambda-1)}{2}\widetilde{[2]}^{2}_{q}L^{2}_{2}a^{2}_{2}]+\frac{\mu(\mu-1)}{2}\lambda^{2}\widetilde{[2]}^{2}_{q}L^{2}_{2}a^{2}_{2}\right\}z^{2}+\ldots,
\end{eqnarray*}
\begin{eqnarray*}
&&\left
\{\frac{2\{\widetilde{\mathcal{D}_{q}}[z\widetilde{\mathcal{D}}_{q}(\mathcal{J}^{\eta}_{q}f)(z)]\}^{\lambda}}{\widetilde{\mathcal{D}}_{q}[\mathcal{J}^{\eta}_{q}f(z)-\mathcal{J}^{\eta}_{q}f(-z)]}\right\}^{1-\mu}=1-\lambda(\mu-1)\widetilde{[2]}^{2}_{q}L_{2}a_{2}z
\notag\\
 &-&\left\{(\mu-1)[(\lambda\widetilde{[3]}_{q}-1)\widetilde{[3]}_{q}L_{3}a_{3}+\frac{\lambda(\lambda-1)}{2}\widetilde{[2]}^{4}_{q}L^{2}_{2}a^{2}_{2}]-\frac{\mu(\mu-1)}{2}\lambda^{2}\widetilde{[2]}^{4}_{q}L^{2}_{2}a^{2}_{2}\right\}z^{2}+\ldots
\end{eqnarray*}
so that
\begin{eqnarray}\label{equ3.7}
&&\left
\{\frac{2z[\widetilde{\mathcal{D}}_{q}(\mathcal{J}^{\eta}_{q}f)(z)]^{\lambda}}{\mathcal{J}^{\eta}_{q}f(z)-\mathcal{J}^{\eta}_{q}f(-z)}\right\}^{\mu}\left
\{\frac{2\{\widetilde{\mathcal{D}_{q}}[z\widetilde{\mathcal{D}}_{q}(\mathcal{J}^{\eta}_{q}f)(z)]\}^{\lambda}}{\widetilde{\mathcal{D}}_{q}[\mathcal{J}^{\eta}_{q}f(z)-\mathcal{J}^{\eta}_{q}f(-z)]}\right\}^{1-\mu}=1+\lambda[\mu-(\mu-1)\widetilde{[2]}_{q}]\widetilde{[2]}_{q}L_{2}a_{2}z\notag\\
 &+&\{(\lambda\widetilde{[3]}_{q}-1)[\mu-(\mu-1)\widetilde{[3]}_{q}]L_{3}a_{3}\notag\\
 &+&\frac{\lambda}{2}\left(\lambda[\mu-(\mu-1)\widetilde{[2]}_{q}]^{2}-[\mu-(\mu-1)\widetilde{[2]}^{2}_{q}]\right)\widetilde{[2]}^{2}_{q}L^{2}_{2}a^{2}_{2}\}z^{2}+\ldots,
\end{eqnarray}
and
\begin{eqnarray*}
&&\left
\{\frac{2w[\widetilde{\mathcal{D}}_{q}(\mathcal{J}^{\eta}_{q}g)(w)]^{\lambda}}{\mathcal{J}^{\eta}_{q}g(w)-\mathcal{J}^{\eta}_{q}g(-w)}\right\}^{\mu}=1-\lambda\mu\widetilde{[2]}_{q}L_{2}a_{2}w+
\notag\\
&&\left\{\mu[(\lambda\widetilde{[3]}_{q}-1)L_{3}(2a^{2}_{2}-a_{3})+\frac{\lambda(\lambda-1)}{2}\widetilde{[2]}^{2}_{q}L^{2}_{2}a^{2}_{2}]+\frac{\mu(\mu-1)}{2}\lambda^{2}\widetilde{[2]}^{2}_{q}L^{2}_{2}a^{2}_{2}\right\}w^{2}+\ldots,
\end{eqnarray*}
\begin{eqnarray*}
&&\left
\{\frac{2\{\widetilde{\mathcal{D}}_{q}[w\widetilde{\mathcal{D}}_{q}(\mathcal{J}^{\eta}_{q}g)(w)]\}^{\lambda}}{\widetilde{\mathcal{D}}_{q}[\mathcal{J}^{\eta}_{q}g(w)-\mathcal{J}^{\eta}_{q}g(-w)]}\right\}^{1-\mu}=1+\lambda(\mu-1)\widetilde{[2]}^{2}_{q}L_{2}a_{2}w
\notag\\
&-&\left\{(\mu-1)[(\lambda\widetilde{[3]}_{q}-1)\widetilde{[3]}_{q}L_{3}(2a^{2}_{2}-a_{3})+\frac{\lambda(\lambda-1)}{2}\widetilde{[2]}^{4}_{q}L^{2}_{2}a^{2}_{2}]-\frac{\mu(\mu-1)}{2}\lambda^{2}\widetilde{[2]}^{4}_{q}L^{2}_{2}a^{2}_{2}\right\}w^{2}+\ldots
\end{eqnarray*}
so that
\begin{eqnarray}\label{equ3.8}
&&\left
\{\frac{2w[\widetilde{\mathcal{D}}_{q}(\mathcal{J}^{\eta}_{q}g)(w)]^{\lambda}}{\mathcal{J}^{\eta}_{q}g(w)-\mathcal{J}^{\eta}_{q}g(-w)}\right\}^{\mu}\left
\{\frac{2\{\widetilde{\mathcal{D}}_{q}[w\widetilde{\mathcal{D}}_{q}(\mathcal{J}^{\eta}_{q}g)(w)]\}^{\lambda}}{\widetilde{\mathcal{D}}_{q}[\mathcal{J}^{\eta}_{q}g(w)-\mathcal{J}^{\eta}_{q}g(-w)]}\right\}^{1-\mu}=1-\lambda[\mu-(\mu-1)\widetilde{[2]}_{q}]\widetilde{[2]}_{q}L_{2}a_{2}w\notag\\
 &+&\{(\lambda\widetilde{[3]}_{q}-1)[\mu-(\mu-1)\widetilde{[3]}_{q}]L_{3}(2a^{2}_{2}-a_{3})+\frac{\lambda}{2}(\lambda[\mu-(\mu-1)\widetilde{[2]}_{q}]^{2}
 \notag\\
 &-&[\mu-(\mu-1)\widetilde{[2]}^{2}_{q}])\widetilde{[2]}^{2}_{q}L^{2}_{2}a^{2}_{2}\}w^{2}+\ldots.
\end{eqnarray}
Therefore, according to (\ref{equ2.5}-\ref{equ2.6})
and(\ref{equ3.5}-\ref{equ3.8}) we have that
\begin{equation}\label{equ3.9}
\lambda[\mu-(\mu-1)\widetilde{[2]}_{q}]\widetilde{[2]}_{q}L_{2}a_{2}=\frac{1}{2}E_{1}c_{1},
\end{equation}
\begin{eqnarray}\label{equ3.10}
(\lambda\widetilde{[3]}_{q}-1)[\mu-(\mu-1)\widetilde{[3]}_{q}]L_{3}a_{3}&+&\frac{\lambda}{2}\left(\lambda[\mu-(\mu-1)\widetilde{[2]}_{q}]^{2}
-[\mu-(\mu-1)\widetilde{[2]}^{2}_{q}]\right)\widetilde{[2]}^{2}_{q}L^{2}_{2}a^{2}_{2}\notag\\
&=&\frac{1}{2}E_{1}\left(c_{2}-\frac{c^{2}_{1}}{2}\right)+\frac{1}{4}E_{2}c^{2}_{1},
\end{eqnarray}
\begin{equation}\label{equ3.11}
-\lambda[\mu-(\mu-1)\widetilde{[2]}_{q}]\widetilde{[2]}_{q}L_{2}a_{2}=\frac{1}{2}E_{1}d_{1}
\end{equation}
and
\begin{eqnarray}\label{equ3.12}
(\lambda\widetilde{[3]}_{q}-1)[\mu-(\mu-1)\widetilde{[3]}_{q}]L_{3}(2a^{2}_{2}-a_{3})&+&\frac{\lambda}{2}\left(\lambda[\mu-(\mu-1)\widetilde{[2]}_{q}]^{2}
-[\mu-(\mu-1)\widetilde{[2]}^{2}_{q}]\right)\widetilde{[2]}^{2}_{q}L^{2}_{2}a^{2}_{2}\notag\\
&=&\frac{1}{2}E_{1}\left(d_{2}-\frac{d^{2}_{1}}{2}\right)+\frac{1}{4}E_{2}d^{2}_{1}.
\end{eqnarray}
From (\ref{equ3.9}) and (\ref{equ3.11}), it infers that
\begin{equation}\label{equ3.13}
a_{2}=\frac{E_{1}c_{1}}{2\lambda[\mu-(\mu-1)\widetilde{[2]}_{q}]\widetilde{[2]}_{q}L_{2}}=-\frac{E_{1}d_{1}}{2\lambda[\mu-(\mu-1)\widetilde{[2]}_{q}]\widetilde{[2]}_{q}L_{2}}
\end{equation}
such that
\begin{equation}\label{equ3.14}
c_{1}=-d_{1}
\end{equation}
and
\begin{equation}\label{equ3.15}
E^{2}_{1}(c^{2}_{1}+d^{2}_{1})=8\lambda^{2}[\mu-(\mu-1)\widetilde{[2]}_{q}]^{2}\widetilde{[2]}^{2}_{q}L^{2}_{2}a^{2}_{2}.
\end{equation}
By (\ref{equ3.10}) and (\ref{equ3.12}), we get that
\begin{eqnarray}\label{equ3.16}
\{2(\lambda\widetilde{[3]}_{q}-1)[\mu-(\mu-1)\widetilde{[3]}_{q}]L_{3}&+&\lambda\left(\lambda[\mu-(\mu-1)\widetilde{[2]}_{q}]^{2}
-[\mu-(\mu-1)\widetilde{[2]}^{2}_{q}]\right)\widetilde{[2]}^{2}_{q}L^{2}_{2}\}a^{2}_{2}\notag\\
&=&\frac{1}{4}(E_{2}-E_{1})(c^{2}_{1}+d^{2}_{1})+\frac{1}{2}E_{1}(c_{2}+d_{2}).
\end{eqnarray}
Therefore, by (\ref{equ3.15}-\ref{equ3.16}) we obtain that
\begin{equation}\label{equ3.17}
a^{2}_{2}=\frac{E^{3}_{1}(c_{2}+d_{2})}{2\Theta},
\end{equation}
where
\begin{eqnarray*}
\Theta:&=&2(\lambda\widetilde{[3]}_{q}-1)[\mu-(\mu-1)\widetilde{[3]}_{q}]E^{2}_{1}L_{3}\\
&+&\lambda\{\lambda[\mu-(\mu-1)\widetilde{[2]}_{q}]^{2}(E^{2}_{1}+2E_{1}-2E_{2})-[\mu-(\mu-1)\widetilde{[2]}^{2}_{q}]E^{2}_{1}\}\widetilde{[2]}^{2}_{q}L^{2}_{2}.
\end{eqnarray*}
Then, it enables us to follow from Lemma \ref{lem1} and (\ref{equ3.15}-\ref{equ3.17})
that
\begin{equation*}
\mid
a_{2}\mid\leq\frac{E_{1}}{\lambda\vert[\mu-(\mu-1)\widetilde{[2]}_{q}]L_{2}\vert\widetilde{[2]}_{q}},
\end{equation*}
\begin{equation*}
\mid a_{2}\mid\leq\sqrt{\frac{2(\vert E_{2}-E_{1}\vert+
E_{1})}{\vert2(\lambda\widetilde{[3]}_{q}-1)[\mu-(\mu-1)\widetilde{[3]}_{q}]L_{3}+\lambda(\lambda[\mu-(\mu-1)\widetilde{[2]}_{q}]^{2}-[\mu-(\mu-1)\widetilde{[2]}^{2}_{q}])\widetilde{[2]}^{2}_{q}L^{2}_{2}\vert}}
\end{equation*}
and
\begin{equation*}
\mid
a_{2}\mid\leq\frac{E_{1}\sqrt{2E_{1}}}{\sqrt{\vert\Theta\vert}}.
\end{equation*}
As a consequence, (\ref{equ3.1}) holds true.

Similarly, from (\ref{equ3.10}), (\ref{equ3.12}) and
(\ref{equ3.14}), it also implies that
\begin{equation}\label{equ3.18}
2(\lambda\widetilde{[3]}_{q}-1)[\mu-(\mu-1)\widetilde{[3]}_{q}]L_{3}(a_{3}-a^{2}_{2})=\frac{1}{2}E_{1}(c_{2}-d_{2}).
\end{equation}
Hence, inserting (\ref{equ3.15}) into (\ref{equ3.18}) we obtain that
\begin{equation}\label{equ3.19}
a_{3}=\frac{E_{1}(c_{2}-d_{2})}{4(\lambda\widetilde{[3]}_{q}-1)[\mu-(\mu-1)\widetilde{[3]}_{q}]L_{3}}
+\frac{E^{2}_{1}(c^{2}_{1}+d^{2}_{1})}{8\lambda^{2}[\mu-(\mu-1)\widetilde{[2]}_{q}]^{2}\widetilde{[2]}^{2}_{q}L^{2}_{2}}.
\end{equation}
Therefore, from Lemma \ref{lem1} it demonstrates that
\begin{equation*}
\mid
a_{3}\mid\leq\frac{E_{1}}{(\lambda\widetilde{[3]}_{q}-1)\vert[\mu-(\mu-1)\widetilde{[3]}_{q}]
L_{3}\vert}
+\frac{E^{2}_{1}}{\lambda^{2}[\mu-(\mu-1)\widetilde{[2]}_{q}]^{2}\widetilde{[2]}^{2}_{q}\vert
L_{2}\vert^{2}}.
\end{equation*}
On the other hand, by (\ref{equ3.16}) and (\ref{equ3.18}) we infer
that
\begin{eqnarray*}
&&a_{3}=\frac{E_{1}(c_{2}-d_{2})}{4(\lambda\widetilde{[3]}_{q}-1)[\mu-(\mu-1)\widetilde{[3]}_{q}]L_{3}}\notag\\
&+&\frac{\frac{1}{4}(E_{2}-E_{1})(c^{2}_{1}+d^{2}_{1})+\frac{1}{2}E_{1}(c_{2}+d_{2})}
{2(\lambda\widetilde{[3]}_{q}-1)[\mu-(\mu-1)\widetilde{[3]}_{q}]L_{3}+\lambda\left(\lambda[\mu-(\mu-1)\widetilde{[2]}_{q}]^{2}
-[\mu-(\mu-1)\widetilde{[2]}^{2}_{q}]\right)\widetilde{[2]}^{2}_{q}L^{2}_{2}}.
\end{eqnarray*}
Thus, from Lemma \ref{lem1} we see that
\begin{eqnarray*}
\vert a_{3}\vert&\leq&\frac{E_{1}}{(\lambda\widetilde{[3]}_{q}-1)\vert[\mu-(\mu-1)\widetilde{[3]}_{q}]L_{3}\vert}\notag\\
&+&\frac{2(\vert E_{2}-E_{1}\vert+E_{1})}
{\vert2(\lambda\widetilde{[3]}_{q}-1)[\mu-(\mu-1)\widetilde{[3]}_{q}]L_{3}+\lambda\left(\lambda[\mu-(\mu-1)\widetilde{[2]}_{q}]^{2}
-[\mu-(\mu-1)\widetilde{[2]}^{2}_{q}]\right)\widetilde{[2]}^{2}_{q}L^{2}_{2}\vert}.
\end{eqnarray*}
\end{proof}

\begin{corollary}\label{cor1} If $f(z)$ given by (\ref{equ1.1}) is in the class
$\widetilde{\mathcal{S^{*}_{\sum}}}^{\eta}_{q}(\lambda;\phi)$,
then
\begin{equation*}
\mid a_{2}\mid\leq\min\left\{\frac{E_{1}}{\lambda\vert
L_{2}\vert\widetilde{[2]}_{q}}, \sqrt{\frac{2(\vert
E_{2}-E_{1}\vert+ E_{1})}{\vert\Gamma\vert}},
\frac{E_{1}\sqrt{2E_{1}}}{\sqrt{\vert\Xi\vert}}\right\}
\end{equation*}
and
\begin{equation*}
\mid a_{3}\mid\leq
\frac{E_{1}}{(\lambda\widetilde{[3]}_{q}-1)\vert
L_{3}\vert}+\min\left\{\frac{E^{2}_{1}}{\lambda^{2}\widetilde{[2]}^{2}_{q}\vert
L_{2}\vert^{2}},
\frac{2(E_{2}-E_{1}\vert+E_{1})}{\vert\Gamma\vert}\right\},
\end{equation*}
where
\begin{equation*}
\Gamma:=2(\lambda\widetilde{[3]}_{q}-1)L_{3}+\lambda(\lambda-1)\widetilde{[2]}^{2}_{q}L^{2}_{2}
\end{equation*}
and
\begin{equation*}
\Xi:=2(\lambda\widetilde{[3]}_{q}-1)E^{2}_{1}L_{3}+\lambda[(\lambda-1)E^{2}_{1}+2\lambda(E_{1}-E_{2})]\widetilde{[2]}^{2}_{q}L^{2}_{2}.
\end{equation*}
\end{corollary}

\begin{corollary}\label{cor2} If $f(z)$ given by (\ref{equ1.1}) is in the class
$\widetilde{\mathcal{C_{\sum}}}^{\eta}_{q}(\mu,\lambda;\phi)$,
then
\begin{equation*}
\mid a_{2}\mid\leq\min\left\{\frac{E_{1}}{\lambda\vert
L_{2}\vert\widetilde{[2]}^{2}_{q}}, \sqrt{\frac{2(\vert
E_{2}-E_{1}\vert+ E_{1})}{\vert\Psi\vert}},
\frac{E_{1}\sqrt{2E_{1}}}{\sqrt{\vert\Upsilon\vert}}\right\}
\end{equation*}
and
\begin{equation*}
\mid a_{3}\mid\leq
\frac{E_{1}}{(\lambda\widetilde{[3]}_{q}-1)\widetilde{[3]}_{q}\vert
L_{3}\vert}+\min\left\{\frac{E^{2}_{1}}{\lambda^{2}\widetilde{[2]}^{4}_{q}\vert
L_{2}\vert^{2}}, \frac{\vert
2(E_{2}-E_{1}\vert+E_{1})}{\vert\Psi\vert}\right\},
\end{equation*}
where
\begin{equation*}
\Psi:=2(\lambda\widetilde{[3]}_{q}-1)\widetilde{[3]}_{q}L_{3}+\lambda(\lambda-1)\widetilde{[2]}^{4}_{q}L^{2}_{2}
\end{equation*}
and
\begin{equation*}
\Upsilon:=2(\lambda\widetilde{[3]}_{q}-1)\widetilde{[3]}_{q}E^{2}_{1}L_{3}+\lambda\left[(\lambda-1)E^{2}_{1}+2\lambda(E_{1}-E_{2})\right]\widetilde{[2]}^{4}_{q}L^{2}_{2}.
\end{equation*}
\end{corollary}

From now on, we pay attention to Fekete-Szeg\"{o} problems for the
class
$\widetilde{\mathcal{S^{*}_{\sum}}}^{\eta}_{q}(\mu,\lambda;\phi)$.

\begin{theorem}\label{thm3}
If $f(z)$ given by (\ref{equ1.1}) is in the class
$\widetilde{\mathcal{S^{*}_{\sum}}}^{\eta}_{q}(\mu,\lambda;\phi)$
and $\varrho\in\mathbb{R}$, then
$$
\mid a_{3}-\varrho a^{2}_{2}\mid\leq\left\{\begin{array}{ll}
\frac{E_{1}}{(\lambda\widetilde{[3]}_{q}-1)\vert[\mu-(\mu-1)\widetilde{[3]}_{q}]
L_{3}\vert},~~&\mbox{if}~~2\vert(1-\varrho)L_{3}[\mu-(\mu-1)\widetilde{[3]}_{q}]\vert
E_{1}^{2}(\lambda\widetilde{[3]}_{q}-1)\leq\vert\Theta\vert,\\
\frac{2\vert1-\varrho\vert E^{3}_{1}}{\vert\Theta\vert},~~
&\mbox{if}~~2\vert(1-\varrho)L_{3}[\mu-(\mu-1)\widetilde{[3]}_{q}]\vert
E_{1}^{2}(\lambda\widetilde{[3]}_{q}-1)\leq\vert\Theta\vert.
\end{array}\right.
$$
where $\Theta$ is the same as in Theorem \ref{thm2}.
\end{theorem}
\begin{proof}[\bf Proof.]
From (\ref{equ3.18}), it infers that
\begin{equation*}
a_{3}-a^{2}_{2}=\frac{E_{1}(c_{2}-d_{2})}{4(\lambda\widetilde{[3]}_{q}-1)[\mu-(\mu-1)\widetilde{[3]}_{q}]L_{3}}.
\end{equation*}
By (\ref{equ3.17}) we obtain that
\begin{eqnarray*}
a_{3}-\varrho
a^{2}_{2}&=&\frac{E_{1}\{2(1-\rho)E_{1}^{2}(\lambda\widetilde{[3]}_{q}-1)[\mu-(\mu-1)\widetilde{[3]}_{q}]L_{3}+\Theta\}c_{2}}{4(\lambda\widetilde{[3]}_{q}-1)[\mu-(\mu-1)\widetilde{[3]}_{q}]
L_{3}\Theta}
\notag\\
&+&\frac{E_{1}\{2(1-\varrho)E_{1}^{2}(\lambda\widetilde{[3]}_{q}-1)[\mu-(\mu-1)\widetilde{[3]}_{q}]L_{3}-\Theta\}d_{2}}{4(\lambda\widetilde{[3]}_{q}-1)[\mu-(\mu-1)\widetilde{[3]}_{q}]
L_{3}\Theta}.
\end{eqnarray*}
Hence,  we know from Lemma \ref{lem1} that
\begin{equation*}
\mid a_{3}-\varrho
a^{2}_{2}\mid\leq\frac{E_{1}}{(\lambda\widetilde{[3]}_{q}-1)\vert[\mu-(\mu-1)\widetilde{[3]}_{q}]
L_{3}\vert}
\end{equation*}
when $ 2\vert(1-\varrho)L_{3}[\mu-(\mu-1)\widetilde{[3]}_{q}]\vert
E_{1}^{2}(\lambda\widetilde{[3]}_{q}-1)\leq\vert\Theta\vert $, or
\begin{equation*}
\mid a_{3}-\varrho a^{2}_{2}\mid\leq\frac{2\vert1-\varrho\vert
E^{2}_{1}}{\vert\Theta\vert}
\end{equation*}
when $2\vert(1-\varrho)L_{3}[\mu-(\mu-1)\widetilde{[3]}_{q}]\vert
E_{1}^{2}(\lambda\widetilde{[3]}_{q}-1)\geq\vert\Theta\vert$.
Then, Theorem \ref{thm3} is completely proved.
\end{proof}

\begin{corollary}\label{cor3}
If $f(z)$ given by (\ref{equ1.1}) is in the class
$\widetilde{\mathcal{S^{*}_{\sum}}}^{\eta}_{q}(\lambda;\phi)$ and
$\varrho\in\mathbb{R}$, then
$$
\mid a_{3}-\varrho a^{2}_{2}\mid\leq\left\{\begin{array}{ll}
\frac{E_{1}}{(\lambda\widetilde{[3]}_{q}-1)
\vert L_{3}\vert},~~&\mbox{if}~~2\vert(1-\varrho)L_{3}\vert E_{1}^{2}(\lambda\widetilde{[3]}_{q}-1)\leq\vert\Xi\vert,\\
\frac{2\vert1-\varrho\vert E^{3}_{1}}{\vert\Xi\vert},~~
&\mbox{if}~~2\vert(1-\varrho)L_{3}\vert
E_{1}^{2}(\lambda\widetilde{[3]}_{q}-1)\geq\vert\Xi\vert,
\end{array}\right.
$$
where
\begin{equation*}
\Xi:=2(\lambda\widetilde{[3]}_{q}-1)E^{2}_{1}L_{3}+\lambda[(\lambda-1)E^{2}_{1}+2\lambda(E_{1}-E_{2})]\widetilde{[2]}^{2}_{q}L^{2}_{2}.
\end{equation*}
\end{corollary}

\begin{corollary}\label{cor4}
If $f(z)$ given by (\ref{equ1.1}) is in the class
$\widetilde{\mathcal{C_{\sum}}}^{\eta}_{q}(\lambda;\phi)$ and
$\varrho\in\mathbb{R}$, then
$$
\mid a_{3}-\varrho a^{2}_{2}\mid\leq\left\{\begin{array}{ll}
\frac{E_{1}}{(\lambda\widetilde{[3]}_{q}-1)\widetilde{[3]}_{q}
\vert L_{3}\vert},~~&\mbox{if}~~2\vert(1-\varrho)L_{3}\vert E_{1}^{2}(\lambda\widetilde{[3]}_{q}-1)\widetilde{[3]}_{q}\leq\vert\Upsilon\vert,\\
\frac{2\vert1-\varrho\vert E^{3}_{1}}{\vert\Upsilon\vert},~~
&\mbox{if}~~2\vert(1-\varrho)L_{3}\vert
E_{1}^{2}(\lambda\widetilde{[3]}_{q}-1)\widetilde{[3]}_{q}\geq\vert\Upsilon\vert,
\end{array}\right.
$$
where
\begin{equation*}
\Upsilon:=2(\lambda\widetilde{[3]}_{q}-1)\widetilde{[3]}_{q}E^{2}_{1}L_{3}
+\lambda[(\lambda-1)E^{2}_{1}+\lambda(E_{1}-E_{2})]\widetilde{[2]}^{4}_{q}L^{2}_{2}.
\end{equation*}
\end{corollary}

\begin{corollary}\label{cor5}
If $f(z)$ given by (\ref{equ1.1}) is in the class
$\widetilde{\mathcal{S^{*}_{\sum}}}^{\eta}_{q}(\mu,\lambda;\phi)$,
then
\begin{equation*}
\mid
a_{3}-a^{2}_{2}\mid\leq\frac{E_{1}}{(\lambda\widetilde{[3]}_{q}-1)\vert[\mu-(\mu-1)\widetilde{[3]}_{q}]
L_{3}\vert}.
\end{equation*}
\end{corollary}

\bigskip

\noindent Pinhong Long

\medskip

\noindent School of Mathematics and Statistics, Ningxia
University, Yinchuan, Ningixa 750021, People's Republic of China

\smallskip

\noindent{\it E-mails:} \texttt{longph@nxu.edu.cn }

\medskip

\noindent Huili Han

\medskip

\noindent School of Mathematics and Computer Science, Ningxia
Normal University, Guyuan, Ningxia 756000, People's Republic of
China

\smallskip

\noindent{\it E-mail:} \texttt{nxhan@126.com }

\medskip

\noindent Halit Orhan

\medskip

\noindent Department of Mathematics, Faculty of Science, Ataturk
University, Erzurum 25240, Turkey

\smallskip

\noindent{\it E-mail:} \texttt{orhanhalit607@gmail.com }

\medskip

\noindent Huo Tang

\medskip

\noindent School of Mathematics and Computer Sciences, Chifeng
University, Chifeng, Inner Mongolia, 024000, People's Republic of
China

\smallskip

\noindent{\it E-mail:} \texttt{thth2009@163.com}

\bigskip

\end{document}